\documentclass[11pt,twoside]{amsart}
\usepackage{latexsym,amssymb,amsmath}
\usepackage[all]{xy}
\usepackage{tikz,pgfplots}
\usepackage{extpfeil}
\usepackage{afterpage}
\usepackage{float}
\usepackage{hyperref}

\makeatletter  
\@namedef{subjclassname@2020}{%
\textup{2020} Mathematics Subject Classification}
\makeatother

\numberwithin{equation}{section}
\newtheorem{theorem}{Theorem}[section]
\newtheorem{lemma}[theorem]{Lemma}
\newtheorem{proposition}[theorem]{Proposition}
\newtheorem{corollary}[theorem]{Corollary}

\theoremstyle{definition}
\newtheorem{definition}[theorem]{Definition} 
\newtheorem{procedure}[theorem]{Procedure} 

\newtheorem{remark}[theorem]{Remark}
\newtheorem{example}[theorem]{Example}

\begin{document}

\title[A duality theorem for the ic-resurgence]
{A duality theorem for the ic-resurgence of edge ideals} 
\thanks{The author was supported by SNI, Mexico.}

\author[R. H. Villarreal]{Rafael H. Villarreal}
\address{
Departamento de
Matem\'aticas\\
Centro de Investigaci\'on y de Estudios
Avanzados del
IPN\\
Apartado Postal
14--740 \\
07000 Mexico City, Mexico
}
\email{vila@math.cinvestav.mx}


\begin{abstract}
The aim of this work is to use linear programming and polyhedral
geometry to prove a duality formula for the
ic-resurgence of edge ideals. We show that the ic-resurgence
of the edge ideal $I$ of a clutter $\mathcal{C}$ and the
ic-resurgence of the edge ideal 
$I^\vee$ of the blocker $\mathcal{C}^\vee$ of $\mathcal{C}$ coincide.
If $\mathcal{C}$ is the clutter of bases of certain uniform matroids, we
recover a formula for the resurgence of $I$, and if $\mathcal{C}$ is a connected
non-bipartite graph with a perfect matching, we show a formula for the
Waldschmidt constant of $I^\vee$.
\end{abstract}

\maketitle 

\section{Introduction}\label{section-intro}
Let $\mathcal C$ be a \textit{clutter} with vertex 
set $V(\mathcal{C})=\{t_1,\ldots,t_s\}$, that is, $\mathcal C$ is a 
family of subsets $E(\mathcal{C})$ of $V(\mathcal{C})$, called edges,
none of which is contained in
another \cite{cornu-book}. We assume that all edges of $\mathcal{C}$ have at least two
vertices. For example, a graph (no multiple edges or loops) is a
clutter.  Regarding each vertex $t_i$ as a
 variable, we consider the polynomial ring $S=K[t_1,\ldots,t_s]$ over
 a field $K$. The monomials of $S$ are denoted
by $t^a:=t_1^{a_1}\cdots t_s^{a_s}$, $a=(a_1,\dots,a_s)$ in $\mathbb{N}^s$, where
$\mathbb{N}=\{0,1,\ldots\}$. The
 \textit{edge ideal} of $\mathcal{C}$, denoted $I(\mathcal{C})$, 
is the ideal of $S$ given by 
$$I(\mathcal{C}):=(\textstyle\{\prod_{t_i\in e}t_i\mid e\in
E(\mathcal{C})\}).$$
\quad The
minimal set of generators of $I(\mathcal{C})$, denoted
$G(I(\mathcal{C}))$, is  the
set of all squarefree monomials $t_e=\prod_{t_i\in e}t_i$ such 
that $e\in E(\mathcal{C})$. 
Any squarefree monomial ideal $I$ of $S$ is the edge
ideal $I(\mathcal{C})$ of a clutter $\mathcal{C}$ with vertex set
$V(\mathcal{C})=\{t_1,\ldots,t_s\}$.   
A set of vertices $C$ of $\mathcal{C}$ is called a \textit{vertex
cover} if every edge of $\mathcal{C}$ contains at least one vertex of
$C$. A  \textit{minimal vertex cover} of $\mathcal{C}$ is a vertex cover which
is minimal with respect to inclusion. The clutter of minimal vertex
covers of $\mathcal{C}$, denoted $\mathcal{C}^\vee$, is called the
\textit{blocker} of $\mathcal{C}$, and the edge ideal 
$I(\mathcal{C}^\vee)$ of $\mathcal{C}^\vee$ is called the
\textit{Alexander dual} of $I(\mathcal{C})$ and is denoted by
$I(\mathcal{C})^\vee$  \cite[p.~221]{monalg-rev}. We assume that 
$|C|\geq 2$ for all $C\in E(\mathcal{C}^\vee)$. 

We denote the edge ideal $I(\mathcal{C})$ of $\mathcal{C}$ by $I$ and 
denote the minimal set of generators of $I$ by
$G(I):=\{t^{v_1},\ldots,t^{v_q}\}$. The \textit{incidence matrix} of $I$, denoted by $A$,
is the $s\times q$ matrix with column vectors $v_1,\ldots,v_q$. This
matrix is the \textit{incidence matrix} of $\mathcal{C}$. The
$\textit{covering polyhedron}$ of $I$, denoted by
$\mathcal{Q}(I)$, is the rational polyhedron
$$
\mathcal{Q}(I):=\{x\mid x\geq 0;\,xA\geq 1\},
$$
where $1=(1,\ldots,1)$. The map
$E(\mathcal{C}^\vee)\rightarrow\{0,1\}^s$, $C\mapsto\sum_{t_i\in C}e_i$, 
induces a bijection between $E(\mathcal{C}^\vee)$ and the set $\{u_1,\ldots,u_m\}$ of integral vertices
of $\mathcal{Q}(I)$ \cite[Corollary~13.1.3]{monalg-rev}. 
Given an integer $n\geq 1$, the $n$-th \textit{symbolic
power} of $I$, denoted $I^{(n)}$, is
given by  \cite[p.~78]{reesclu}: 
\begin{equation}\label{jun21-21-1}
I^{(n)}=(\{t^a\mid a/n\in\mathcal{Q}(I^\vee)\})=(\{t^a\mid
\langle a, u_i\rangle\geq n\mbox{ for }i=1,\ldots,m\}),
\end{equation}
where $\langle\ ,\, \rangle$ denotes the 
ordinary inner product. 
The \textit{Newton
polyhedron} of $I$, denoted ${\rm NP}(I)$, is the 
integral polyhedron 
\begin{equation*}
{\rm NP}(I)=\mathbb{R}_+^s+{\rm 
conv}(v_1,\ldots,v_q),
\end{equation*}
where $\mathbb{R}_+=\{\lambda\in\mathbb{R}\mid \lambda\geq
0\}$. It is well known that ${\rm NP}(I)$ is equal to 
$$\mathcal{Q}(\overline{B}):=\{x\mid 
x\geq 0;\,x\overline{B}\geq 1\},
$$
where $\overline{B}$ is the rational matrix whose columns are
precisely the vertices $u_1,\ldots,u_p$ of
the covering polyhedron $\mathcal{Q}(I)$ of $I$
\cite[Proposition~3.5(b)]{reesclu}. 
Note that $m\leq p$, that
is, not all vertices of $\mathcal{Q}(I)$ are integral. The \textit{integral closure} of
 $I^n$, denoted $\overline{I^n}$, can be described as 
\begin{equation}\label{jun21-21}
\overline{I^n}=
(\{t^a\mid a/n\in{\rm NP}(I)\})=(\{t^a\mid
\langle a, u_i\rangle\geq n\mbox{ for }i=1,\ldots,p\})
\end{equation}
for every $n\geq 1$  \cite[Proposition~3.5(a)]{reesclu}. 
If $I^n=\overline{I^n}$ for all $n\geq 1$, $I$ is said to be
\textit{normal}. 

The resurgence and asymptotic resurgence of ideals were introduced 
in \cite{resurgence,asymptotic-resurgence}. The resurgence of an ideal
relative to the integral closure was introduced in
\cite{Francisco-TAMS}. We define the \textit{resurgence},
\textit{asymptotic resurgence}, and \textit{ic-resurgence} of $I$ to be 
\begin{align*}
&\rho(I):=\sup \left\{{n}/{r}
\ \left. \right|\, 
I^{(n)}\not\subset I^r\right\},\\ 
&\ \ \ \ \ \widehat{\rho}(I):=\sup
\left \{ {n}/{r}\ 
\left.\right|\, 
I^{(nt)}\not\subset I^{rt} \text{ for all } t\gg 0\right\},\\ 
&\ \ \ \ \ \ \ \ \ \ \ \rho_{ic}(I):=\sup
\left \{ {n}/{r}\ 
\left.\right|\, 
I^{(n)}\not\subset\overline{I^{r}}\right\},\
\text{respectively}.
\end{align*}
\quad In general, ${\rho}(I)\geq {\rho}_{ic}(I)$ and $\widehat{\rho}(I)={\rho}_{ic}(I)$ 
\cite[Corollary~4.14]{Francisco-TAMS}. In particular
$\widehat{\rho}(I)$ is ${\rho}(I)$ if $I$ is normal. 
A related well studied invariant is $\widehat{\alpha}(I)$, the \textit{Waldschmidt constant} of
$I$ \cite{Cooper-symbolic,Francisco-TAMS}:  
\begin{equation*}
\widehat{\alpha}(I):=\lim_{n\rightarrow\infty}\frac{\alpha(I^{(n)})}{n},
\end{equation*}
where $\alpha(I^{(n)})$ is the 
least degree of a minimal generator
of $I^{(n)}$. Let $B$ be the
incidence matrix of $I^\vee$. Note that the columns of $B$ are
$u_1,\ldots,u_m$. The number $\widehat{\alpha}(I)$ is the
optimal value of the linear program
\begin{enumerate}
\item[] {\rm minimize} $y_1+\cdots+y_s$

\item[] {\rm subject to}

\item[] $yB\geq 1$ {\rm and }$y\geq 0$
\end{enumerate}
which is attained at a vertex of
$\mathcal{Q}(I^\vee)$ \cite[Theorem~3.2]{Waldschmidt-Bocci-etal}. 
If $\alpha(\mathcal{Q}(I))$ is the minimum of all
$|a|=\sum_{i=1}^sa_i$ with $a=(a_1,\ldots,a_s)$ a vertex of
$\mathcal{Q}(I)$, then
$\widehat{\alpha}(I^\vee)=\alpha(\mathcal{Q}(I))$ 
\cite[Corollary~6.3]{Cooper-symbolic}. 

It is known that the ic-resurgence
of the edge ideal $I(G)$ of a graph $G$ can be expressed in terms of the 
vertices of $\mathcal{Q}(I(G)^\vee)$. 
The following formula was given in \cite[Theorem~3.12]{Francisco-TAMS}:
\begin{equation*}
\rho_{ic}(I(G))=\widehat{\rho}(I(G))=\frac{2}{\widehat{\alpha}(I(G))}=
\frac{2}{\alpha(\mathcal{Q}(I(G)^\vee))}\geq
1.
\end{equation*}
\quad Our main  result is a formula for the ic-resurgence
$\rho_{ic}(I)$ of $I$ in terms of
the vertices of $\mathcal{Q}(I)$ and $\mathcal{Q}(I^\vee)$
(Theorem~\ref{duality-ic-resurgence}). 

As we now explain computing the ic-resurgence $\rho_{ic}(I)$ of the
ideal $I$ is a linear-fractional programming 
problem \cite{boyd}. Let 
$V(\mathcal{Q}(I))=\{u_1,\ldots,u_p\}$ be the vertex
set of $\mathcal{Q}(I)$, we can write $u_i=\gamma_i/d_i$, 
$\gamma_i\in\mathbb{N}^s\setminus\{0\}$, $d_i\in\mathbb{N}_+$ for all
$i$ (Section~\ref{duality-section}).  
Recall that we may assume that $d_i=1$ for $i=1,\ldots,m$, i.e.,
$u_1,\ldots,u_m$ are the characteristic vectors 
of the minimal vertex covers of $\mathcal{C}$ and $m\leq p$. 
By Eqs.~\eqref{jun21-21-1}-\eqref{jun21-21}, a monomial
$t^a$ is in $I^{(n)}\setminus\overline{I^r}$ if and only if
$a/n\in\mathcal{Q}(I^\vee)$ and $a/r\notin{\rm NP}(I)$, that is,  $t^a$ is in
$I^{(n)}\setminus\overline{I^r}$  if and only if $a=(a_1,\ldots,a_s)$ satisfies
\begin{equation}\label{mar14-20}
\langle a,u_i\rangle\geq n\mbox{ for }i=1,\ldots,m\mbox{ and }\langle
a,\gamma_j\rangle\leq rd_j-1\mbox{ for some }1\leq j\leq p.
\end{equation}
\quad Let $x_1,\ldots,x_s$ be variables that correspond to
the entries of $a$ and let $x_{s+1},x_{s+2}$ be two extra variables
that correspond to $n$ and $r$, respectively. By Eq.~\eqref{mar14-20}, 
$t^a$ is in $I^{(n)}\setminus\overline{I^r}$ if and only if there
exists $1\leq j\leq p$ such that $(a,n,r)$ is a 
feasible point for the linear-fractional program:
\begin{align}
&\text{maximize }\ \ h_j(x)=\frac{x_{s+1}}{x_{s+2}}&&&&\nonumber\\
&\text{subject to }\ \langle(x_1,\ldots,x_s),u_i\rangle-x_{s+1}\geq
0,\ i=1,\ldots,m,\, x_{s+1}\geq 1 &&&&\label{lfp-ic-resurgence}\\
&\quad \quad\quad\quad\quad x_i\geq 0,\, i=1,\ldots,s&&&&\nonumber \\
&\quad \quad\quad\quad\quad  d_jx_{s+2}-\langle(x_1,\ldots,x_s),\gamma_j\rangle
\geq 1,\, x_{s+2}\geq 1&&&&\nonumber 
\end{align}

This type of program can be solved using linear
programming \cite[Section~4.3.2,
p.~151]{boyd}. The linear-fractional program of
Eq.~\eqref{lfp-ic-resurgence} is equivalent to the 
linear program of Eq.~\eqref{lp-ic-resurgence} below
\cite[p.~18]{intclos}. The next result
solves the problem of computing the ic-resurgence of $I$. An
algorithm to compute $\rho_{ic}(I)$ has been implemented in 
\textit{Macaulay}$2$ using its interface to 
\textit{Normaliz} \cite{normaliz2} (see \cite[Procedure~A.3,
Algorithm~A.4]{intclos}). 
We use the next
result to prove a duality
formula for the ic-resurgence of $I$ that implies the equality 
$\rho_{ic}(I)=\rho_{ic}(I^\vee)$.

\begin{theorem}\label{lp-resurgence-formula}
\cite[Theorem~5.3]{intclos}
For each $1\leq j\leq p$, let $\rho_j$
be the optimal value of the 
following linear
program with variables $y_1,\ldots,y_{s+3}$. Then,
$\rho_{ic}(I)=\max\{\rho_j\}_{j=1}^p$ 
and $\rho_j$ is 
attained at a rational vertex of the polyhedron $\mathcal{P}_j$ of feasible points of
Eq.~\eqref{lp-ic-resurgence}. 
\begin{align}
&\text{maximize }\ \ g_j(y)=y_{s+1}&&&&\nonumber\\
\ \ &\text{subject to }\ \langle(y_1,\ldots,y_s),u_i\rangle-y_{s+1}\geq
0,\ i=1,\ldots,m,\, y_{s+1}\geq y_{s+3} &&&&\label{lp-ic-resurgence}\\
&\quad \quad\quad\quad\quad y_i\geq 0,\, i=1,\ldots,s,\, y_{s+3}\geq
0&&&&\nonumber \\
&\quad \quad\quad\quad\quad  d_jy_{s+2}-\langle(y_1,\ldots,y_s),\gamma_j\rangle
\geq y_{s+3},\, y_{s+2}=1.&&&&\nonumber 
\end{align}
\end{theorem}

First, we determine all vertices $y$ of $\mathcal{P}_j$ with
$y_{s+1}>0$ and relate them with the vertices of $\mathcal{Q}(I)$ and $\mathcal{Q}(I^\vee)$.
To introduce our result, let $V(\mathcal{Q}(I^\vee))=\{v_1,\ldots,v_{p_1}\}$ be the vertex
set of $\mathcal{Q}(I^\vee)$, we can write $v_i=\delta_i/f_i$,
$\delta_i\in\mathbb{N}^s\setminus\{0\}$, 
$f_i\in\mathbb{N}_+$ for all $i$ (Section~\ref{duality-section}).
Recall that we may assume that $f_i=1$ for $i=1,\ldots,q$, i.e.,
$v_1,\ldots,v_q$ are the characteristic vectors 
of the edges of $\mathcal{C}$ and $q\leq p_1$. We introduce the Rees
cones of $I$ and $I^\vee$ as a device to compute the vertices of
$\mathcal{Q}(I)$ and $\mathcal{Q}(I^\vee)$, and the linear constraints that
define the Newton polyhedron of $I$
(Theorems~\eqref{irred-rep-rc}-\eqref{reesclu-theorem-d}). 
Then, we collect some
basic properties of $\mathcal{P}_j$ showing that $\mathcal{P}_j$
could be unbounded, and proving that $\rho_j\leq d_j$ 
and $\langle\gamma_i,\delta_\ell\rangle\geq d_i\geq 1$ for all
$i,\ell$ (Lemma~\ref{well-defined-invariant}). 

We come to our main auxiliary result.

\noindent \textbf{Theorem~\ref{vertices-pj}.}\textit{ Let $y=(y_1,\ldots,y_{s+3})$ be a point
in $\mathbb{R}^{s+3}$ with $y_{s+1}>0$ and let $1\leq j\leq p$, then $y$ is a vertex of
$\mathcal{P}_j$ if and only if $y$ has one of the following two forms:
\begin{enumerate}
\item[(a)] $y=(y_{s+1}(\delta_k/f_k),\, y_{s+1},\, 1,\, 0)$, where 
$y_{s+1}=d_jf_k/\langle
\gamma_j,\delta_k\rangle$ and $1\leq k\leq p_1$, or 
\item[(b)] $y=(y_{s+1}(\delta_k/f_k),\, y_{s+1},\, 1,\, y_{s+1})$, where 
$y_{s+1}=d_jf_k/(\langle
\gamma_j,\delta_k\rangle+f_k)$ and $1\leq k\leq p_1$.
\end{enumerate}
}

The first evidence that the equality $\rho_{ic}(I)=\rho_{ic}(I^\vee)$
could be true came from \cite{Francisco-TAMS} and
\cite{intclos,Jayanthan-Kumar-Mukundan} 
where it is shown that for the edge ideal $I(G)$ of a perfect graph
$G$, one has       
$$
\rho_{ic}(I(G))=\frac{2(\omega(G)-1)}{\omega(G)}=\rho_{ic}(I(G)^\vee),
$$
respectively, where $\omega(G)$ is the clique number of $G$, that is,
$\omega(G)$ is the number of vertices in a maximum complete subgraph
of $G$. If $G$ is a perfect graph, then $\rho(I(G)^\vee)$ is
equal to $\rho_{ic}(I(G)^\vee)$ because $I(G)^\vee$ is normal 
\cite[Theorem~2.10]{perfect}. 

We come to our main result.

\noindent \textbf{Theorem~\ref{duality-ic-resurgence}.} (Duality
formula)\textit{ 
If $I$ is the edge ideal of a clutter $\mathcal{C}$, then 
$$
\frac{1}{\rho_{ic}(I)}=\min\left\{\langle u,v\rangle\mid u\in V(\mathcal{Q}(I)),\,
v\in V(\mathcal{Q}(I^\vee))\right\},
$$
where $V(\mathcal{Q}(I))$ is the vertex
set of $\mathcal{Q}(I)$. In particular, $\rho_{ic}(I)=\rho_{ic}(I^\vee)$.
}

A covering polyhedron is
\textit{integral} if and only if it has only integral vertices
\cite[p.~232]{Schr}. 
The ic-resurgence of $I$ classifies the integrality of 
the covering polyhedron $\mathcal{Q}(I)$ because $\rho_{ic}(I)=1$ if and only if
$\mathcal{Q}(I)$ is integral (Proposition~\ref{feb5-22}). As a
consequence of the duality formula for the ic-resurgence, we recover
the fact that $\mathcal{Q}(I)$ is integral if and only if
$\mathcal{Q}(I^\vee)$ is integral \cite[Theorem~1.17]{cornu-book} 
(Corollary~\ref{feb8-22-1}), and recover two results of Jayanthan, 
Kumar and Mukundan \cite[Theorems~4.8 and
5.3]{Jayanthan-Kumar-Mukundan}  
showing that when $I(G)$ is the edge ideal of a graph 
$G$ the following conditions are equivalent 
\begin{align*}
&{\rm (a)}\ \rho(I(G))=1;\quad{\rm
(b)}\ \rho_{ic}(I(G))=1; 
\quad {\rm (c)}\ G\mbox{ is bipartite};\quad{\rm (d)}\
\rho(I(G)^\vee)=1; 
\end{align*}
see Corollary~\ref{feb9-22}. Then, we show a formula for the ic-resurgence
of the sum of two edge ideals of clutters generated by monomials in disjoint sets of
variables (Corollary~\ref{feb8-22}). A similar formula is known for
the  Waldschmidt constant (Remark~\ref{similar-w}). If $G$ is a connected
non-bipartite graph with a perfect matching, we show that $\widehat{
\alpha}(I(G)^\vee)=|V(G)|/2$ (Proposition~\ref{w-formula-dual}).

To compute the vertices of the covering polyhedron $\mathcal{Q}(I)$ of
the edge ideal $I$ of a clutter $\mathcal{C}$, we have used
Theorems~\ref{irred-rep-rc} and \ref{reesclu-theorem}, see
Procedure~\ref{ic-resurgence-example1.m2}. Another way to compute the
vertices of $\mathcal{Q}(I)$ is to find the extreme rays, i.e., the
$1$-dimensional faces of the following polyhedral cone
\begin{equation*}
{\rm SC}(I^\vee):=\{x\in\mathbb{R}^{s+1}\vert\, x\geq 0;\, \langle
x,(v_i,-1)\rangle\geq 0\ \mbox{ for } \ i=1,\ldots,q\}, 
\end{equation*}
see \cite[Proposition 3.15]{intclos}. 
The cone ${\rm SC}(I^\vee)$ is called the \textit{Simis cone}
of $I^\vee$ \cite{normali}. For a discussion on how to find all vertices of a general
polyhedron we refer to \cite{avis-fukuda,avis-fukuda1}.

The algebraic properties and invariants of the ideal $I=I_{d,s}$, generated by
all squarefree monomials of $S$ of degree $d$ in $s$ variables, were studied in
\cite{adrian-ainv,Waldschmidt-Bocci-etal,Geramita-Harbourne,Vi4,matrof}. This ideal is
normal \cite{Vi4} and its Alexander dual $I^\vee$ is equal to $I_{s-d+1,s}$. The minimal generators
of $I$ correspond to the bases of a matroid of rank $d$. 
The statement about resurgence in Theorem~\ref{sq-Veronese} is 
due originally to Lampa-Baczy\'nska and Malara
\cite[Theorem~C]{Lampa-Malara}. In a more general form, it appears in
a paper by Geramita, Harbourne, Migliore
and Nagel \cite[Theorem~4.8]{Geramita-Harbourne}. The
Waldschmidt constant of $I$ was computed in \cite[Theorem~7.5]{Waldschmidt-Bocci-etal}.
 Using the Rees cones of
$I$ and $I^\vee$, we compute the vertices 
of $\mathcal{Q}(I)$ and $\mathcal{Q}(I^\vee)$ and use
Theorem~\ref{duality-ic-resurgence} to recover the formulas for $\rho(I)$ and
$\widehat{\alpha}(I)$. 

\noindent
\textbf{Theorem~\ref{sq-Veronese}.}\ \cite{Waldschmidt-Bocci-etal,Lampa-Malara}\textit{
The resurgence and ic-resurgence of $I_{s,d}$ and $I_{s,d}^\vee$ are
given by
$$
\rho(I_{d,s})=\rho_{ic}(I_{d,s})=\rho_{ic}((I_{d,s})^\vee)=
\frac{d(s-d+1)}{s}=\frac{d}{\widehat{\alpha}(I_{d,s})}=
\frac{s-d+1}{\widehat{\alpha}((I_{d,s})^\vee)}.
$$
}

In Section~\ref{examples-section} we present examples illustrating 
our results. Then in Appendix~\ref{Appendix} we give the procedures
for \textit{Macaulay}$2$
\cite{mac2}, that are used in the 
examples. 

For unexplained
terminology and additional information,  we refer to 
\cite{mc8,huneke-swanson-book,Vas1,bookthree} for the theory of
integral closure, 
\cite{Herzog-Hibi-book,monalg-rev} for the theory of
edge ideals and monomial ideals, and \cite{korte,Schr,Schr2} for
combinatorial 
optimization and integer
programing.  

\section{Preliminaries}\label{section-prelim} 
In this section we introduce a few results from
polyhedral geometry and commutative algebra. We continue to employ 
the notations and definitions used in Section~\ref{section-intro}. 

Given $a\in {\mathbb R}^s\setminus\{0\}$  and 
$c\in {\mathbb R}$, the \textit{affine hyperplane} $H{(a,c)}$ 
and the \textit{positive closed halfspace} $H^+{(a,c)}$ 
bounded by $H{(a,c)}$ are defined as
\[
H{(a,c)}:=\{x\in{\mathbb R}^s\vert\, \langle x,a\rangle=c\}\ \mbox{ and
}\ H^+{(a,c)}:=\{x\in{\mathbb R}^s\vert\, 
\langle x,a\rangle\geq c\}.
\]
\quad If $c=0$, $H_a$ will denote $H{(a,c)}$ and $H_a^+$ will denote
$H^+{(a,c)}$. 
If $a$ and $c$ are rational, $H^+{(a,c)}$ 
is called a \textit{rational closed halfspace}. 

A \textit{rational polyhedron}
is a subset of ${\mathbb R}^s$ which is the 
intersection of a finite number of rational closed halfspaces of 
$\mathbb{R}^s$. Let $\Gamma$ be a subset of $\mathbb{R}^s$. The cone 
generated by $\Gamma$, 
denoted ${\mathbb R}_+\Gamma$, is the set of all 
linear combinations of $\Gamma$ with coefficients in the set
$\mathbb{R}_+$ of nonnegative real numbers.
The \textit{finite basis theorem} asserts that a subset $\mathcal{U}$ of $\mathbb{R}^s$ is a 
rational polyhedron if and only if
$$\mathcal{U}=\mathbb{R}_+\Gamma+\mathcal{P},$$  
where $\mathbb{R}_+\Gamma$ is a cone generated by a finite
set $\Gamma$ of rational points and $\mathcal{P}$ is the convex hull ${\rm conv}(\mathcal{A})$ of a
finite set $\mathcal{A}$ of rational
points \cite[Corollary~7.1b]{Schr}. The
computer programs \textit{Normaliz}
\cite{normaliz2} and  \textit{PORTA} \cite{porta} can be used to switch between these 
two representations.

\begin{theorem}\cite[Theorem~4.3]{Seceleanu-packing}\label{jun12-21-1}
Let $I$ be a squarefree monomial ideal. If $\mathcal{Q}(I)$ is 
integral, then $\widehat{\alpha}(I)=\alpha(I)$ and
$\widehat{\alpha}(I^\vee)=\alpha(I^\vee)$.
\end{theorem}

\begin{theorem}\label{ntf-char} Let
$I$ be a squarefree monomial ideal. The following hold. 
\begin{enumerate}
\item[(a)] \cite{clutters} 
$I^n=I^{(n)}$ for all $n\geq 1$ if and only if $\mathcal{Q}(I)$ is integral
and $I$ is normal. 
\item[(b)] \cite{clutters,Trung} $\overline{I^n}=I^{(n)}$
for all $n\geq 1$ if and only if $\mathcal{Q}(I)$ is integral.
\end{enumerate}
\end{theorem}

\section{Duality formula for the ic-resurgence}\label{duality-section}
In this part we give a duality formula for the ic-resurgence of 
edge ideals of clutters. 
To avoid repetitions, we continue to employ 
the notations and definitions used in Sections~\ref{section-intro} 
and \ref{section-prelim}. 

Let $S=K[t_1,\ldots,t_s]$ be a polynomial ring over a field
$K$, let $I\subset S$ be the edge ideal of a clutter
$\mathcal{C}$, and let
$G(I)=\{t^{v_1},\ldots,t^{v_q}\}$ be the minimal set of generators of
$I$. To study $\rho_{ic}(I)$,  
we need a convenient way to determine the polyhedron
$\mathcal{P}_j$ of feasible points of the linear program of
Theorem~\ref{lp-resurgence-formula} for $1\leq j\leq p$. 
Our approach is based on the computation of the supporting hyperplanes
of the \textit{Rees cone} of $I$ which is the
finitely generated rational cone defined as \cite{normali}: 
\begin{equation}\label{rees-cone-eq}
{\rm RC}(I):=\mathbb{R}_+\{e_1,\ldots,e_s,(v_1,1),\ldots,(v_q,1)\}.
\end{equation}

\begin{theorem}\cite[Proposition~1.1.51,
Theorem~14.1.1]{monalg-rev}\label{irred-rep-rc} The Rees cone of
$I$ has a unique irreducible representation 
\begin{equation}\label{supp-hyp}
{\rm RC}(I)=\left(\bigcap_{i=1}^{s+1}H^+_{e_i}\right)\bigcap
\left(\bigcap_{i=1}^mH^+_{(\gamma_i,-d_i)}\right)
\bigcap\left(\bigcap_{i=m}^p H^+_{(\gamma_i,-d_i)}\right),
\end{equation}
where none of the closed halfspaces can be omitted from the
intersection, $d_i=1$ for $i=1,\ldots,m$, $t^{\gamma_1},\ldots,t^{\gamma_m}$ are the 
minimal generators of $I^\vee$, $\gamma_i\in\mathbb{N}^s\setminus\{0\}$ for $i>m$, 
$d_i\in\mathbb{N}\setminus\{0,1\}$ for $i>m$, and the non-zero entries of
$(\gamma_i,-d_i)$ are relatively prime for all $i$. 
\end{theorem}

The hyperplanes defining the closed halfspaces of Eq.~\eqref{supp-hyp} are the
\textit{supporting hyperplanes} of the Rees cone of $I$, and the $\gamma_i$'s
and $d_i$'s can be computed using \textit{Normaliz} \cite{normaliz2}. 
Thus, for each $1\leq j\leq p$, we can determine the
polyhedron $\mathcal{P}_j$.

The following theorem justifies the a priori
naming coincidence between the constants appearing in
Theorem~\ref{irred-rep-rc}
 and the  description given for the $u_i$ in the introduction before Eq.~\eqref{mar14-20}.

\begin{theorem}\cite[Theorem~3.1]{reesclu}\label{reesclu-theorem} 
The vertex set of $\mathcal{Q}(I)$ is
$V(\mathcal{Q}(I))=\{\gamma_1/d_1,\ldots,\gamma_p/d_p\}$.
\end{theorem}

The last two results say that finding the supporting hyperplanes
of ${\rm RC}(I)$ is equivalent to finding the vertices of
$\mathcal{Q}(I)$. 

For use below we
set $u_i=\gamma_i/d_i=\gamma_i$ for $i=1,\ldots,m$. 
To study the vertices of $\mathcal{P}_j$ we need the following dual
versions of Theorems~\ref{irred-rep-rc} and \ref{reesclu-theorem}.
The Rees cone of $I^\vee$ is given by 
\begin{equation}\label{rees-cone-eq-dual}
{\rm RC}(I^\vee)=\mathbb{R}_+\{e_1,\ldots,e_s,(u_1,1),\ldots,(u_m,1)\}.
\end{equation}

\begin{theorem}\label{irred-rep-rcd}
The Rees cone of $I^\vee$ has a unique irreducible representation 
\begin{equation}\label{supp-hyp-dual}
{\rm RC}(I^\vee)=\left(\bigcap_{i=1}^{s+1}H^+_{e_i}\right)\bigcap
\left(\bigcap_{i=1}^qH^+_{(\delta_i,-f_i)}\right)
\bigcap\left(\bigcap_{i=q+1}^{p_1} H^+_{(\delta_i,-f_i)}\right),
\end{equation}
where $f_i=1$ for $i=1,\ldots,q$,
$t^{\delta_1},\ldots,t^{\delta_q}$ are the 
minimal generators of $I$, $\delta_i\in\mathbb{N}^s\setminus\{0\}$
for $i>q$, 
$f_i\in\mathbb{N}\setminus\{0,1\}$ for $i>q$, and the non-zero entries of
$(\delta_i,-f_i)$ are relatively prime for all $i$. 
\end{theorem}

\begin{theorem}\label{reesclu-theorem-d} 
The vertex set of $\mathcal{Q}(I^\vee)$ is
$V(\mathcal{Q}(I^\vee))=\{\delta_1/f_1,\ldots,\delta_{p_1}/f_{p_1}\}$.
\end{theorem}

In what follows we set $v_i=\delta_i/f_i$ for $i=1,\ldots,p_1$ and 
$u_i=\gamma_i/d_i$ for $i=1,\ldots,p$. Recall 
that $\{t^{v_1},\ldots,t^{v_q}\}$ is the minimal generating set of the
edge ideal $I=I(\mathcal{C})$ of $\mathcal{C}$.
The minimal generating set
for $I^\vee$ is $\{t^{u_1},\ldots,t^{u_m}\}$ and
$I^\vee=I(\mathcal{C}^\vee)$ is the edge
ideal of the blocker $\mathcal{C}^\vee$ of $\mathcal{C}$.

\begin{lemma}\label{well-defined-invariant} 
Let $I$ be the edge ideal of a clutter $\mathcal{C}$. The following
hold. 
\begin{enumerate}
\item[(a)] Given $\gamma_i$ and $\delta_\ell$, 
there are $u\in\{u_i\}_{i=1}^m$ and $v\in\{v_i\}_{i=1}^q$
such that $\gamma_i\geq u$ and $\delta_\ell\geq v$ componentwise.
\item[(b)] $\langle\gamma_i,\delta_\ell\rangle\geq d_i\geq 1$ and 
$\langle\gamma_i,\delta_\ell\rangle\geq f_\ell\geq 1$ for all
$i,\ell$. 
\item[(c)] If $\rho_j$ is  the optimal value of the linear program of
Eq.~\eqref{lp-ic-resurgence}, then $\rho_j\leq d_j$.
\item[(d)] If $\gamma_j=(\gamma_{j,1},\ldots,\gamma_{j,s})$ and
$\gamma_{j,k}=0$ for some $k$, then $\mathcal{P}_j$ is an unbounded
polyhedron.
\item[(e)] If $\gamma_j=(\gamma_{j,1},\ldots,\gamma_{j,s})$ and
$\gamma_{j,k}>0$ for $k=1,\ldots,s$, then $\mathcal{P}_j$ is a bounded
polyhedron.
\end{enumerate}
\end{lemma}

\begin{proof} (a)-(b) Given $a\in\mathbb{N}^s$,
$a=(a_1,\ldots,a_s)$, we set
$F_a={\rm supp}(t^a)=\{t_k\mid a_k>0\}$. From Eq.~\eqref{supp-hyp},
$\langle v_n,\gamma_i\rangle\geq d_i$ for $n=1,\ldots,q$. Then,
$F_{v_n}\bigcap F_{\gamma_i}\neq\emptyset$ for $n=1,\ldots,q$. As 
$F_{v_1},\ldots,F_{v_q}$ are the edges of $\mathcal{C}$,
$F_{\gamma_i}$ is a vertex cover of $\mathcal{C}$, and consequently 
$F_{\gamma_i}$ contains a minimal vertex cover $C$ of $\mathcal{C}$.
Similarly, from Eq.~\eqref{supp-hyp-dual}, $\langle u_n,\delta_\ell\rangle\geq
f_\ell$ for $n=1,\ldots,m$. Then,
$F_{u_n}\bigcap F_{\delta_\ell}\neq\emptyset$ for $n=1,\ldots,m$. As 
$F_{u_1},\ldots,F_{u_m}$ are the edges of $\mathcal{C}^\vee$,
$F_{\delta_\ell}$ is a vertex cover of $\mathcal{C}^\vee$, and consequently 
$F_{\delta_\ell}$ contains a minimal vertex cover $D$ of
$\mathcal{C}^\vee$. Note that $D$ is an edge of $\mathcal{C}$ because
$(\mathcal{C}^\vee)^\vee=\mathcal{C}$. Let $u=\sum_{t_i\in C}e_i$
and $v=\sum_{t_i\in D}e_i$ be the characteristic vectors of $C$ and $D$, respectively.
Then, $\delta_\ell\geq v$ and $\gamma_i\geq u$. Hence, using 
Eq.~\eqref{supp-hyp} and Eq.~\eqref{supp-hyp-dual}, we obtain
\begin{align*}
&\langle\gamma_i,\delta_\ell\rangle\geq\langle \gamma_i,v\rangle\geq
d_i\geq 1\mbox{ and }\langle\gamma_i,\delta_\ell\rangle\geq\langle u, \delta_\ell\rangle\geq
f_\ell\geq 1.
\end{align*}

(c) Let
$y=(y_1,\ldots,y_{s+3})$ 
be any point in $\mathcal{P}_j$, that is, $y$ is feasible for 
Eq.~\eqref{lp-ic-resurgence}. 
By part (a), one can pick $u\in\{u_i\}_{i=1}^m$ such that
$\gamma_j\geq u$. Then, using the constraints that define $\mathcal{P}_j$, we
get 
$$y_{s+1}\leq\langle
(y_1,\ldots,y_s),u\rangle\leq\langle(y_1,\ldots
y_s),\gamma_j\rangle\leq d_j-y_{s+3}\leq d_j.
$$
\quad Thus, $y_{s+1}\leq d_j$, and consequently $\rho_j\leq d_j$.

(d) For simplicity, we may assume that $k=1$. Setting
$x=(x_1,\ldots,x_{s+3})=\ell e_1+e_{s+2}$ for $\ell\geq 0$ and using
that $\gamma_{j,1}=0$, one obtains
\begin{align*}
&\langle (x_1,\ldots,x_s), u_i\rangle\geq x_{s+1}=0,\ i=1,\ldots,m,\
x_{s+1}\geq x_{s+3}\\
&x_i\geq 0,\ i=1,\ldots,s,\ x_{s+3}\geq 0\\
& d_jx_{s+2}-\langle(x_1,\ldots,x_s),\gamma_j\rangle
=d_j-x_1\gamma_{j,1}=d_j\geq x_{s+3}=0,\, x_{s+2}=1,
\end{align*}
and consequently $x\in\mathcal{P}_j$ for all $\ell\geq 0$. This proves
that $\mathcal{P}_j$ is unbounded.

(e) Let
$y=(y_1,\ldots,y_{s+3})$ 
be any point in $\mathcal{P}_j$. It suffices to show that $y_i\leq
d_j$ for $i=1,\ldots,s+3$. By part (a), one can pick $u\in\{u_i\}_{i=1}^m$ such that
$\gamma_j\geq u$. As $\gamma_{j}$ is integral,
$\gamma_{j,k}\geq 1$ for all $k$. Then, using the constraints that
define $\mathcal{P}_j$, we get  
\begin{align*}
&y_{i}\leq\langle(y_1,\ldots
y_s),\gamma_j\rangle\leq d_j-y_{s+3}\leq d_j,\ i=1,\ldots,s\\
&y_{s+3}\leq y_{s+1}\leq\langle
(y_1,\ldots,y_s),u\rangle\leq\langle(y_1,\ldots
y_s),\gamma_j\rangle\leq d_j-y_{s+3}\leq d_j,
\end{align*}
and the proof is complete.
\end{proof}

\begin{theorem}\label{vertices-pj} Let $\mathcal{P}_j$ be the polyhedron of feasible points of
Eq.~\eqref{lp-ic-resurgence}, $1\leq j\leq p$, and let $y=(y_1,\ldots,y_{s+3})$ be a point
in $\mathbb{R}^{s+3}$ with $y_{s+1}>0$, then $y$ is a vertex of
$\mathcal{P}_j$ if and only if $y$ has one of the following two forms:
\begin{enumerate}
\item[(a)] $y=(y_{s+1}(\delta_k/f_k),\, y_{s+1},\, 1,\, 0)$, where 
$y_{s+1}=d_jf_k/\langle
\gamma_j,\delta_k\rangle$ and $1\leq k\leq p_1$, or 
\item[(b)] $y=(y_{s+1}(\delta_k/f_k),\, y_{s+1},\, 1,\, y_{s+1})$, where 
$y_{s+1}=d_jf_k/(\langle
\gamma_j,\delta_k\rangle+f_k)$ and $1\leq k\leq p_1$.
\end{enumerate}
\end{theorem}

\begin{proof} $\Rightarrow$) Assume that $y=(y_1,\ldots,y_{s+3})$ is
a vertex of $\mathcal{P}_j$ with $y_{s+1}>0$. The constraints that
define $\mathcal{P}_j$ can be written as
\begin{align}
&\langle y,(-u_i,1,0,0)\rangle\leq
0,\ i=1,\ldots,m,\, \langle y,-e_{s+1}+e_{s+3}\rangle\leq 0,\nonumber\\
&\langle y,-e_i \rangle\leq 0,\, i=1,\ldots,s,\, \langle y,-e_{s+3}
\rangle\leq 0,\label{constraints-vectors}\\
&\langle y, (\gamma_j,0,0,0)-d_je_{s+2}+e_{s+3}\rangle\leq 0,\,
\langle y, e_{s+2}\rangle\leq 1,\, \langle y, -e_{s+2}\rangle\leq -1.\nonumber 
\end{align}
\quad As $y$ is a vertex of $\mathcal{P}_j$, by \cite[Corollary~1.1.47]{monalg-rev},
there are $s+3$ linearly independent constraints that are satisfied
with equality. The hyperplane
$H_{(y_1,\ldots,y_{s+1})}$ of $\mathbb{R}^{s+1}$ contains at most $s$
linearly independent vectors because $y_{s+1}>0$. Therefore, there are exactly $s$ constraints of the form  
\begin{align}
&\langle y,(-u_i,1,0,0)\rangle\leq
0,\ i=1,\ldots,m,\label{cqj-1} \\ 
&\langle y,-e_i \rangle\leq 0,\, 
i=1,\ldots,s,\label{cqj-0}
\end{align}
that are satisfied with equality because the remaining constraints
that define $\mathcal{P}_j$ are
\begin{align}
&\langle y,-e_{s+1}+e_{s+3}\rangle\leq 0,\label{cqj1}\\
&\langle y,-e_{s+3}
\rangle\leq 0,\label{cqj2}\\
&\langle y, (\gamma_j,0,0,0)-d_je_{s+2}+e_{s+3}\rangle\leq
0,\label{cqj3}\\ 
&\langle y, e_{s+2}\rangle\leq 1,\, \langle y, -e_{s+2}\rangle\leq
-1,\label{cqj4} 
\end{align}
and the constraints of Eqs.~\eqref{cqj1} and \eqref{cqj2} cannot
hold  simultaneously with equality. Then, recalling that $y_{s+2}=1$
in $\mathcal{P}_j$, one has  
\begin{equation}\label{cqj5}
\langle(y_1,\ldots,y_s),\gamma_j\rangle=d_j-y_{s+3}
\end{equation}
because the constraint of Eq.~\eqref{cqj3} must hold with equality,
and either $y_{s+1}=y_{s+3}$ if Eq.\eqref{cqj1} holds 
with equality or $y_{s+3}=0$ if Eq.\eqref{cqj2} holds with equality. Let 
$$
\mathcal{B}=\{(-w_i,1)\}_{i=1}^{r-1}\textstyle\bigcup\{(-w_i,0)\}_{i=r}^s
$$ 
be the set of vectors in $\mathbb{R}^{s+1}$ that correspond to the
$s$ linearly independent constraints of Eqs.~\eqref{cqj-1} and
\eqref{cqj-0} that are satisfied with equality. Note that 
$\mathcal{B}$ is linearly independent if and only if 
$\mathcal{B}'=\{(w_i,1)\}_{i=1}^{r-1}\bigcup\{(w_i,0)\}_{i=r}^s$ is linearly
independent. Consider the hyperplane 
$$
H=H_{(y_1,\ldots,y_s,-y_{s+1})}.
$$ 
\quad The set $\mathcal{B}'$ is
contained in $H$, and $\mathcal{B}'$ is also contained in the Rees
cone ${\rm RC}(I^\vee)$ of $I^\vee$ because $w_i\in\{u_1,\ldots,u_m\}$
for $1\leq i\leq r-1$ and $w_i\in\{e_1,\ldots,e_s\}$ for $r\leq i\leq
s$. Then, $F=H\cap{\rm RC}(I^\vee)$ is a facet of
the Rees cone ${\rm RC}(I^\vee)$ of $I^\vee$. Hence, by
Theorem~\ref{irred-rep-rcd} and \cite[Theorem~3.2.1]{webster}, we
obtain that $F=H_{(\delta_k,-f_k)}\cap {\rm RC}(I^\vee)$ for some
$1\leq k\leq p_1$, and consequently
$$(\delta_k,-f_k)=\lambda(y_1,\ldots,y_s,-y_{s+1})$$
for some $\lambda\in\mathbb{R}$. Thus,
$(y_1,\ldots,y_s)/y_{s+1}=\delta_k/f_k$ for some $1\leq k\leq p_1$.

Case (I) $y_{s+3}=0$. Then, dividing Eq.~\eqref{cqj5} by
$y_{s+1}$, we get 
$$
\frac{\langle
\delta_k,\gamma_j\rangle}{f_k}=\frac{d_j}{y_{s+1}},
$$
and solving for $y_{s+1}$ gives $y_{s+1}=d_jf_k/\langle
\gamma_j,\delta_k\rangle$. Thus, $y$ is as in (a).

Case (II) $y_{s+1}=y_{s+3}$. Then, dividing Eq.~\eqref{cqj5} by
$y_{s+1}$, we get 
$$
\frac{\langle \delta_k,\gamma_j\rangle+f_k}{f_k}=\frac{\langle
\delta_k,\gamma_j\rangle}{f_k}+1=\frac{d_j}{y_{s+1}},
$$
and solving for $y_{s+1}$ gives $y_{s+1}=d_jf_k/(\langle
\gamma_j,\delta_k\rangle+f_k)$. Thus, $y$ is as in (b).

$\Leftarrow$) Assume that $y$ is as in (a). As
$(y_1,\ldots,y_s)/y_{s+1}=\delta_k/f_k$ is a vertex of
$\mathcal{Q}(I^\vee)$, it follows that $y$ satisfies 
the constraints of Eqs.~\eqref{cqj-1}-\eqref{cqj-0} and that 
$s$ of them occur with equality and are linearly independent. 
Let $\Gamma$ be the set of vectors in $\mathbb{R}^{s+3}$ that
correspond to these independent constraints. Then,
$\langle\alpha,y\rangle=0$ for all $\alpha\in\Gamma$. The following
three constraints are also satisfied with equality
\begin{align}
&\langle y,-e_{s+3}
\rangle\leq 0,\nonumber\\
&\langle y, (\gamma_j,0,0,0)-d_je_{s+2}+e_{s+3}\rangle\leq
0,\nonumber\\ 
&\langle y, e_{s+2}\rangle\leq 1.
\nonumber 
\end{align}
\quad As $y\in\mathcal{P}_j$, to prove that $y$ is a vertex of $\mathcal{P}_j$ we need only show that the set 
$$
\Gamma\textstyle\cup\{-e_{s+3},\,
(\gamma_j,0,0,0)-d_je_{s+2}+e_{s+3},\, e_{s+2}\}
$$
is linearly independent or equivalently that 
$\Gamma\textstyle\cup\{e_{s+3},(\gamma_j,0,0,0),e_{s+2}\}$ is
linearly independent. This follows noticing that
$\Gamma\cup\{e_{s+2},e_{s+3}\}$ is linearly independent and observing that $(\gamma_j,0,0,0)$ cannot be
in $\mathbb{R}(\Gamma\cup\{e_{s+2},e_{s+3}\})$, the linear span of
$\Gamma\cup\{e_{s+2},e_{s+3}\}$, because $\langle
y,(\gamma_j,0,0,0)\rangle=d_j$ and $\langle\alpha,y\rangle=0$ for all
$\alpha\in\Gamma$. 

Assume that $y$ is as in (b). As
$(y_1,\ldots,y_s)/y_{s+1}=\delta_k/f_k$ is a vertex of
$\mathcal{Q}(I^\vee)$, it follows that $y$ satisfies 
the constraints of Eqs.~\eqref{cqj-1}-\eqref{cqj-0} and that 
$s$ of them occur with equality and are linearly independent. 
Let $\Gamma$ be the set of vectors in $\mathbb{R}^{s+3}$ that
correspond to these independent constraints. Then,
$\langle\alpha,y\rangle=0$ for all $\alpha\in\Gamma$. The following
constraints are satisfied with equality
\begin{align}
&\langle y,-e_{s+1}+e_{s+3}
\rangle\leq 0,\nonumber\\
&\langle y, (\gamma_j,0,0,0)-d_je_{s+2}+e_{s+3}\rangle\leq
0,\nonumber\\ 
&\langle y, e_{s+2}\rangle\leq 1.
\nonumber 
\end{align}
\quad As $y\in\mathcal{P}_j$, to prove that $y$ is a vertex of $\mathcal{P}_j$ we need only show that the set 
$$
\Gamma\textstyle\cup\{-e_{s+1}+e_{s+3},(\gamma_j,0,0,0)-d_je_{s+2}+e_{s+3},e_{s+2}\}
$$
is linearly independent or equivalently that the set 
$$
\Gamma\textstyle\cup\{-e_{s+1}+e_{s+3},\, (\gamma_j,0,0,0)+e_{s+3},\, e_{s+2}\}
$$ 
is linearly independent. Clearly the set 
$\Gamma\cup\{-e_{s+1}+e_{s+3},e_{s+2}\}$ is linearly independent. 
Hence, it suffices to notice that $(\gamma_j,0,0,0)$ cannot be 
in $\mathbb{R}(\Gamma\cup\{-e_{s+1}+e_{s+3},e_{s+2}\})$ because $\langle
y,(\gamma_j,0,0,0)\rangle$ is equal to $d_j-y_{s+1}$, $\langle\alpha,y\rangle=0$ for all
$\alpha\in\Gamma$, and $d_j-y_{s+1}\neq 0$ (Lemma~\ref{well-defined-invariant}). 
\end{proof}

\begin{theorem}{\rm(Duality formula)}\label{duality-ic-resurgence}
If $I=I(\mathcal{C})$ is the edge ideal of a clutter $\mathcal{C}$
and $\rho_{ic}(I)$ is the ic-resurgence of $I$, then
$$
\frac{1}{\rho_{ic}(I)}=\min\left\{\langle u,v\rangle\mid u\in V(\mathcal{Q}(I)),\,
v\in V(\mathcal{Q}(I^\vee))\right\},
$$
where $V(\mathcal{Q}(I))$ is the vertex
set of $\mathcal{Q}(I)$. In particular, $\rho_{ic}(I)=\rho_{ic}(I^\vee)$.
\end{theorem}

\begin{proof} The sets of vertices of $\mathcal{Q}(I)$ and
$\mathcal{Q}(I^\vee)$ are $\{\gamma_i/d_i\}_{i=1}^p$ and
$\{\delta_\ell/f_\ell\}_{\ell=1}^{p_1}$, respectively. This follows 
from Theorems~\ref{reesclu-theorem} and \ref{reesclu-theorem-d}. 
Then, noticing the equalities 
$$
1/\min\{\langle\gamma_i/d_i,\delta_\ell/f_\ell\rangle\}_{i,\ell}=
\max\{1/\langle\gamma_i/d_i,\delta_\ell/f_\ell\rangle\}_{i,\ell}=
\max\left\{(d_if_\ell)/\langle\gamma_i,\delta_\ell\rangle\right\}_{i,\ell},
$$
it suffices to prove the following equality 
\begin{equation}\label{jan30-22}
\rho_{ic}(I)=
\max\left\{(d_if_\ell)/\langle\gamma_i,\delta_\ell\rangle\right\}_{i,\ell}.
\end{equation}
\quad  Given integers $1\leq i\leq p$ and $1\leq \ell\leq p_1$, by
Theorem~\ref{vertices-pj}, there is a vertex $y_{i,\ell}$ of
$\mathcal{P}_i$ such that the $(s+1)$-th entry $(y_{i,\ell})_{s+1}$ of
$y_{i,\ell}$ is equal to
$(d_if_\ell)/\langle\gamma_i,\delta_\ell\rangle$. As $y_{i,\ell}$ is 
in $\mathcal{P}_i$, we obtain that 
$\rho_i\geq(d_if_\ell)/\langle\gamma_i,\delta_\ell\rangle$, and 
consequently one has 
\begin{align}\label{jan30-22-1}
&\rho_{ic}(I)\geq\rho_i\geq d_if_\ell/\langle
\gamma_i,\delta_\ell\rangle\ \forall\ i,\ell.
\end{align}
\quad By Theorems~\ref{lp-resurgence-formula} and \ref{vertices-pj}, 
$\rho_{ic}(I)=\rho_j$
for some $1\leq j\leq p$ and there is $y=(y_1,\ldots,y_{s+3})$ a
vertex of $\mathcal{P}_j$ such that 
$y_{s+1}=\rho_{ic}(I)$, and either 
$$y_{s+1}=d_jf_k/\langle
\gamma_j,\delta_k\rangle\ or\ y_{s+1}=d_jf_k/(\langle
\gamma_j,\delta_k\rangle+f_k)
$$ 
for some $k$. The second equality cannot occur because 
$$d_jf_k/(\langle
\gamma_j,\delta_k\rangle+f_k)<d_jf_k/\langle
\gamma_j,\delta_k\rangle\leq\rho_{ic}(I).
$$
\quad Thus, one has $y_{s+1}=d_jf_k/\langle
\gamma_j,\delta_k\rangle$ and, by Eq.~\eqref{jan30-22-1}, we get 
equality in Eq.~\eqref{jan30-22}.
\end{proof}

\begin{proposition}\label{feb5-22} If $I$ is a squarefree monomial ideal, then 
$\rho_{ic}(I)\geq 1$ with equality if and only if
$\mathcal{Q}(I)$ is integral. 
\end{proposition}

\begin{proof} By \cite[Corollaries~4.14 and 4.16]{Francisco-TAMS},
$\rho_{ic}(I)\geq 1$ with equality if and only if
$I^{(n)}=\overline{I^n}$ for every $n\geq 1$. Therefore, by
Theorem~\ref{ntf-char}, $\rho_{ic}(I)=1$ if and only if
$\mathcal{Q}(I)$ is integral. 
\end{proof}

\begin{corollary}\cite[Theorem~1.17]{cornu-book}\label{feb8-22-1} 
$\mathcal{Q}(I)$ is integral if and only if
$\mathcal{Q}(I^\vee)$ is integral.
\end{corollary}
\begin{proof} If $\mathcal{Q}(I)$ is integral, by
Proposition~\ref{feb5-22}, $\rho_{ic}(I)=1$. Then, by 
Theorem~\ref{duality-ic-resurgence}, $\rho_{ic}(I^\vee)=1$. Thus, 
by Proposition~\ref{feb5-22}, $\mathcal{Q}(I^\vee)$ is integral. 
The converse follows by replacing $I$ with $I^\vee$, in the first part
of the proof, and recalling that
$(I^\vee)^\vee$ is equal to $I$. 
\end{proof}

\begin{corollary}\cite[Theorems~4.8 and
5.3]{Jayanthan-Kumar-Mukundan}\label{feb9-22} Let $I=I(G)$ be the
edge ideal of a  
graph. The following conditions are equivalent.
\begin{align*}
&{\rm (a)}\ \rho(I)=1;\quad{\rm
(b)}\ \rho_{ic}(I)=1;
\quad {\rm (c)}\ G\mbox{ is bipartite};\quad{\rm (d)}\ \rho(I^\vee)=1. 
\end{align*}
\end{corollary}

\begin{proof} 
(a)$\Rightarrow$(b) Assume that $\rho(I)=1$. Then, 
$\rho(I)\geq\rho_{ic}(I)\geq 1$, and $\rho_{ic}(I)=1$.

(b)$\Rightarrow$(c) Assume that $\rho_{ic}(I)=1$. Then, by Proposition~\ref{feb5-22}, we get
that $\mathcal{Q}(I)$ is integral. Hence, by
\cite[Proposition~14.3.39]{monalg-rev}, $G$ is bipartite. 

(c)$\Rightarrow$(d) Assume that $G$ is bipartite. Then, by
\cite[Corollary~2.6]{alexdual} and
\cite[Proposition~14.3.6]{monalg-rev}, $I^\vee$ is normal and
$\mathcal{Q}(I^\vee)$ integral. Hence, by Proposition~\ref{feb5-22}, we get 
$\rho(I^\vee)=\rho_{ic}(I^\vee)=1$.

(d)$\Rightarrow$(a) Assume that $\rho(I^\vee)=1$. Then, 
$\rho(I^\vee)\geq\rho_{ic}(I^\vee)\geq 1$, and $\rho_{ic}(I^\vee)=1$.
Then, by the duality formula and Proposition~\ref{feb5-22}, we get
that $\mathcal{Q}(I)$ is integral and $\rho_{ic}(I)=1$. Then, by
\cite[Proposition~14.3.39]{monalg-rev} and
\cite[Theorem~5.9]{ITG}, $G$ is bipartite and $I$ is
normal. Hence, $\rho(I)=\rho_{ic}(I)=1$. 
\end{proof}

\begin{corollary}\label{feb8-22} If $I_1$ and $I_2$ are squarefree monomial ideals of
$S$ generated by monomials in disjoint sets of variables, then 
$$ 
\rho_{ic}(I_1+I_2)=\max\{\rho_{ic}(I_1),\, \rho_{ic}(I_2)\}.
$$
\end{corollary}

\begin{proof} We may assume that $\mathcal{C}_1$ and $\mathcal{C}_2$
are clutters with vertex sets
$\{t_1,\ldots,t_r\}$ and $\{t_{r+1},\ldots,t_s\}$ such that
$I_i$ is the edge ideal of $\mathcal{C}_i$ for $i=1,2$. If 
$\mathcal{C}=\mathcal{C}_1\cup\mathcal{C}_2$ and $I=I_1+I_2$, then
$I$ is the edge ideal of $\mathcal{C}$ and $C$ is a minimal vertex
cover of $\mathcal{C}$ if and only if $C=C_1\cup C_2$ with $C_i$ a
minimal vertex cover of $\mathcal{C}_i$ for $i=1,2$. Hence,
$I^\vee=I_1^\vee I_2^\vee$ and, by \cite[Proposition~3.5]{Jayanthan-Kumar-Mukundan}, 
$\rho_{ic}(I^\vee)$ is equal to
$\max\{\rho_{ic}(I_1^\vee),\rho_{ic}(I_2^\vee)\}$. Therefore, using
Theorem~\ref{duality-ic-resurgence}, we get
$$
\rho_{ic}(I)=\rho_{ic}(I^\vee)=\rho_{ic}(I_1^\vee
I_2^\vee)=\max\{\rho_{ic}(I_1^\vee),\rho_{ic}(I_2^\vee)\}=\max\{\rho_{ic}(I_1),\rho_{ic}(I_2)\}.
$$
\quad Thus, $\rho_{ic}(I)=\max\{\rho_{ic}(I_1),\rho_{ic}(I_2)\}$ and the
proof is complete.
\end{proof}

\begin{remark}\label{similar-w} There is a similar formula for 
$\widehat{\alpha}(I_1+I_2)$ \cite[Corollary~7.10]{Waldschmidt-Bocci-etal}: 
$$ 
\widehat{\alpha}(I_1+I_2)=\min\{\widehat{\alpha}(I_1),\,
\widehat{\alpha}(I_2)\}.
$$
\end{remark}

\begin{definition}\cite[p.~541]{monalg-rev} Let $a=(a_i)\neq 0$ be a 
vector in $\mathbb{N}^s$ and let $b\in\mathbb{N}$. If $a, b$ satisfy 
$\langle a,v_i\rangle\geq b$ for $i=1,\ldots,q$, we say that 
$a$ is a $b$-{\it cover\/} of
$\mathcal{C}$. 
\end{definition}

The notion of a $b$-cover 
occurs in combinatorial optimization \cite[Chapter~77,
p.~1378]{Schr2} and algebraic 
combinatorics \cite{normali,cover-algebras}. 

\begin{definition}\rm A $b$-cover $a$ of $\mathcal{C}$ is
called {\it reducible\/} if there exists an $i$-cover $c$ and a $j$-cover
$d$ of $\mathcal{C}$ 
such that $a=c+d$ and $b=i+j$. If $a$ is not
reducible, we call $a$ {\it irreducible\/}.    
\end{definition}

\begin{lemma}\cite[Lemma~1.8]{covers}\label{aug25-07} If $(\gamma_k,-d_k)$ is any
of the vectors of 
Eq.~\eqref{supp-hyp}, then $\gamma_k$ is an irreducible $d_k$-cover of
$\mathcal{C}$. 
\end{lemma}

\begin{definition}\label{konig-def}
Let $G$ be a graph. A set $P$ of pairwise disjoint edges of $G$ is a
\textit{perfect matching} of $G$ if
$V(G)=\bigcup_ {e\in P}e$. A set of vertices of $G$ is 
\textit{stable} if no two of them are adjacent.
\end{definition} 

\begin{proposition}\label{w-formula-dual} Let $G$ be a connected non-bipartite graph with a
perfect matching and let $I(G)$ be its edge ideal. Then,
$\widehat{\alpha}(I(G)^\vee)=|V(G)|/2$.
\end{proposition}

\begin{proof} Let $u_i=\gamma_i/d_i$ be any vertex of $\mathcal{Q}(I)$. 
By Theorem~\ref{reesclu-theorem}, $(\gamma_i,-d_i)$
occurs in Eq.~\eqref{supp-hyp}. Hence, by Lemma~\ref{aug25-07}, 
$\gamma_i$ is an irreducible $d_i$-cover of
$G$. Therefore, using the classification of
the irreducible covers of a graph \cite[Theorem~1.7]{covers}, 
$\gamma_i$ and $d_i$ have one of the following forms
\begin{enumerate}
\item[(a)] $d_i=1$ and $\gamma_i=\sum_{t_j\in C}e_j$ for
some $C\in E(G^\vee)$, that is, $1\leq i\leq m$,
\item[(b)] $d_i=2$ and $\gamma_i=(1,\ldots,1)$,
\item[(c)] $d_i=2$ and up to permutation of vertices 
\begin{equation}\label{mar17-22}
\gamma_i=(\underbrace{0,\ldots,0}_{|\mathfrak{A}|},
\underbrace{2,\ldots,2}_{|N_G(\mathfrak{A})|},
{1,\ldots,1})
\end{equation}
for some stable set of vertices $\mathfrak{A}\neq\emptyset$ of
$G$, where $N_G(\mathfrak{A})$ is the neighbor set of $\mathfrak{A}$.
\end{enumerate}
\quad The incidence matrix of $I$ has rank $s=|V(G)|$ because $G$ is
connected and non-bipartite \cite[Lemma~2.1]{join}. Hence, 
$\beta=(1/2,\ldots,1/2)$ is a vertex of $\mathcal{Q}(I)$. Indeed, note
that $\beta A\geq 1$, where $A$ is the incidence matrix of $I$. There
are $s$ linearly independent columns of $A$, say $v_1,\ldots,v_s$,
such that $\langle\beta,v_i\rangle=1$ for $i=1,\ldots,s$. Thus,
$\beta$ is a basic feasible solution of the system $xA\geq 1;\,x\geq
0$ and, by \cite[Corollary~1.1.49]{monalg-rev}, $\beta$ is a vertex of
$\mathcal{Q}(I)$. Then, by \cite[Theorem~3.2]{Waldschmidt-Bocci-etal}, one has
$$\widehat{\alpha}(I^\vee)=\min\{|u_k|\colon u_k\in V(\mathcal{Q}(I))\}
\leq |\beta|=s/2.
$$ 
\quad Therefore, it suffices to show that $|u_i|\geq s/2$. If
$u_i=\gamma_i/d_i$ is as in (a), then $|u_i|\geq 
s/2$ because $G$ has a perfect matching and any minimal vertex cover $C$ of $G$ has at least $s/2$
vertices of $G$. If $u_i$ is as in (c), by Eq.~\eqref{mar17-22}, we obtain
$$
|u_i|=|\gamma_i/2|=|N_G(\mathfrak{A})|+(1/2)(s-|\mathfrak{A}|-|N_G(\mathfrak{A})|)=
(1/2)|N_G(\mathfrak{A})| +(s/2)-(1/2)|\mathfrak{A}|.
$$
\quad As the graph $G$ has a perfect matching and the set $\mathfrak{A}$ contains no
edges of $G$ since $\mathfrak{A}$ is stable, one has 
$|N_G(\mathfrak{A})|\geq|\mathfrak{A}|$. Thus, $|u_i|\geq s/2$. 
\end{proof}

\begin{lemma}\label{intersection-number}
Let $d,k,\ell,s$ be integers such that $1\leq d\leq \ell\leq
s$ and $s-d+1\leq k\leq s$. If $A$ and $B$ are subsets of
$\{1,\ldots,s\}$, $k=|A|$, $\ell=|B|$, then
\begin{equation}\label{feb27-22-1}
\frac{|A\cap B|}{(k-s+d)(\ell-d+1)}\geq\frac{s}{d(s-d+1)}
\end{equation}
with equality if $A=B=\{1,\ldots,s\}$.
\end{lemma}

\begin{proof} We can write $\ell=d+i$ and $|A\cup B|=s-\epsilon$, where
$0\leq i\leq s-d$ and $\epsilon\geq 0$. As $|A\cap B|$ is equal to $k+\ell-|A\cup
B|$, we obtain that Eq.~\eqref{feb27-22-1} is equivalent to  
\begin{equation*}
(k+d+i-s+\epsilon)d(s-d+1)\geq 
s(k-s+d)(i+1).
\end{equation*}
\quad By factoring out $i$, this inequality is equivalent to 
\begin{equation}\label{feb27-22-3}
\epsilon d(s-d+1)+(k+d-s)(d-1)(s-d)\geq
i\left[s(k-s+d)-d(s-d+1)\right].
\end{equation}
\quad We set $f=s(k-s+d)-d(s-d+1)$. If $f\leq 0$, the inequality of
Eq.~\eqref{feb27-22-3} holds because the left-hand side is
non-negative and the right-hand side is $if$. Thus, we may assume that $f\geq 1$. Since $s-d\geq i$,
one has $(s-d)f\geq if$, and we need only prove the inequality
\begin{equation}\label{feb27-22-4}
\epsilon d(s-d+1)+(k+d-s)(d-1)(s-d)\geq
(s-d)\left[s(k-s+d)-d(s-d+1)\right].
\end{equation}
\quad This inequality is equivalent to
$$
d(s-d+1)(\epsilon+(s-d))\geq(s-d)(s-(d-1))(k-s+d).
$$
\quad By cancelling out $s-d+1$, the proof reduces to showing that 
$$
d(\epsilon+(s-d))\geq(s-d)(k-s+d). 
$$
\quad As $d(s-d)$ appears as a summand on both sides of this inequality it suffices to
notice that $d\epsilon\geq(s-d)(k-s)$ because $k\leq s$. 
\end{proof}

\begin{theorem}\cite{Waldschmidt-Bocci-etal,Lampa-Malara}\label{sq-Veronese}
If $I_{d,s}$ is the $d$-th squarefree Veronese ideal of $S$ generated
by all squarefree monomials of $S$ of degree $d$ in $s$ variables, then 
\begin{equation}\label{mar18-22}
\rho(I_{d,s})=\rho_{ic}(I_{d,s})=\rho_{ic}((I_{d,s})^\vee)=\frac{d(s-d+1)}{s}=\frac{d}{\widehat{\alpha}(I_{d,s})}=
\frac{s-d+1}{\widehat{\alpha}((I_{d,s})^\vee)}.
\end{equation}
\end{theorem}

\begin{proof} We set $I=I_{d,s}$ and $J=I_{s-d+1,s}$. The primary
decomposition of $I$ is 
\begin{equation}\label{feb28-22-1}
I=\bigcap_{1\leq j_1<\cdots<j_{s-d+1}\leq
s}(t_{j_1},\ldots,t_{j_{s-d+1}}).
\end{equation}
\quad Then, $(I_{d,s})^\vee=I_{s-d+1,s}$ and $J=(I_{d,s})^\vee$. 
If $d=s$, then $I=(t_1\cdots t_s)$, $I^\vee=(t_1,\ldots,t_s)$, the
vertex set of $\mathcal{Q}(I)$ is $\{e_1,\ldots,e_s\}$, and the vertex
set of $\mathcal{Q}(I^\vee)$ is $\{e_1+\cdots+e_s\}$. Then, by
\cite[Theorem~3.2]{Waldschmidt-Bocci-etal} and
Theorem~\ref{duality-ic-resurgence}, we get
$\widehat{\alpha}(I^\vee)=1$, $\widehat{\alpha}(I)=s$, and
$\rho_{ic}(I)=\rho_{ic}(I^\vee)=1$. Then, Eq.~\eqref{mar18-22} holds.
For similar reasons Eq.~\eqref{mar18-22} holds if $d=1$. 

Thus, we may assume that $2\leq d<s$. We
claim that the
vertex sets of $\mathcal{Q}(I)$ and $\mathcal{Q}(I^\vee)$ are
\begin{align}
&V(Q(I))=\left\{\left.\frac{e_{i_1}+\cdots+e_{i_k}}{k-s+d}\right| 
s-d+1\leq k\leq s,\ 1\leq i_1<\cdots<i_k\leq
s\right\},\label{eq1-mar1}\\
&V(Q(I^\vee))=\left\{\left.\frac{e_{j_1}+\cdots+e_{j_{\ell}}}{\ell-d+1}\right| 
d\leq \ell\leq s,\ 1\leq j_1<\cdots<j_\ell\leq
s\right\}\label{eq2-mar1}.
\end{align}
\quad To prove the first equality, we use the Rees cone ${\rm RC}(I)$
of $I$ defined in Eq.~\eqref{rees-cone-eq}. The
second equality will follow by using the Rees cone ${\rm RC}(J)$ of $J$
or by replacing $I_{d,s}$ with $I_{s-d+1,s}$. Next we prove that the irreducible
representation of ${\rm RC}(I)$ is given by
\begin{equation}\label{supp-hyp-matroid}
{\rm
RC}(I)=\left(\bigcap_{i=1}^{s+1}H^+_{e_i}\right)\bigcap
\Bigg(
\bigcap_{\begin{array}{c}
\scriptstyle 1 \leq i_1<\cdots<i_{k}\leq
s \vspace{-1mm}\\ 
\scriptstyle s-d+1\leq k\leq s\end{array}
}H^+_{(e_{i_1}+\cdots+e_{i_k},-(k-s+d))}
\Bigg).
\end{equation}
\quad First we show that the non-trivial closed halfspaces on the
right-hand side
of this equality occur in the irreducible representation of ${\rm
RC}(I)$. Consider the hyperplane
$$
H=H_{(e_{i_1}+\cdots+e_{i_k},-(k-s+d))},
$$
where $1 \leq i_1<\cdots<i_{k}\leq
s$ and $s-d+1\leq k\leq s$. To show that $H^+$ occurs
in the irreducible representation of ${\rm RC}(I)$ we need only 
show that the set $F=H\cap{\rm RC}(I)$ is a facet of ${\rm RC}(I)$
\cite[Theorem~3.2.1]{webster}.  
Recall that ${\rm RC}(I)$ is the cone $\mathbb{R}_+\mathcal{A}'$
generated by the set 
$$
\mathcal{A}'=\{e_1,\ldots,e_s\}\cup\{
e_{j_1}+\cdots+e_{j_d}+e_{s+1}\mid 1\leq j_1<\cdots<j_d\leq s
\}. 
$$
\quad That $H$ is a supporting hyperplane of the Rees cone, that is, ${\rm
RC}(I)\subset H^+$, follows using that for any sequence of integers $1\leq
j_1<\cdots<j_d\leq s$ one has
\begin{align*}
&\langle
e_{i_1}+\cdots+e_{i_k},e_{j_1}+\cdots+e_{j_d}\rangle=|\{i_1,\ldots,i_k\}\cap\{j_1,\ldots,j_d\}|=
\nonumber\\
&k+d-|\{i_1,\ldots,i_k\}\cup\{j_1,\ldots,j_d\}|\geq
k+d-s.\nonumber
\end{align*}
\quad Hence, to prove that $F$ is a facet of 
${\rm RC}(I)$, it suffices to prove that $F$ contains $s$ linearly
independent vectors in $\mathcal{A}'$. Recall that $d-s+k<k$ because
$d<s$. Take any subset $A_1$ of
$\{i_1,\ldots,i_k\}$ with $d-s+k$ elements and any subset $A_2$ of
$\{1,\ldots,s\}\setminus\{i_1,\ldots,i_k\}$ with $s-k$ elements. Then
\begin{align}
&\big\langle\sum_{i\in A_1\cup A_2}e_i,\, e_{i_1}+\cdots+e_{i_k}\big \rangle=|(A_1\cup
A_2)\cap\{i_1,\ldots,i_k\}|=d-s+k,\mbox{ and }\label{feb28-22}\\
&e_i\in H \mbox{ for }i\in \{1,\ldots,s\}\setminus\{i_1,\ldots,i_k\}.
\label{feb28-22-bis}
\end{align}
\quad Let $I_{d-s+k, k}$ be the $(d-s+k)$-th squarefree Veronese
ideal in the $k$ variables $\{t_{i_1},\ldots,t_{i_k}\}$. The incidence matrix of
$I_{d-s+k, k}$ has rank $k$ because $d-s+k<k$ \cite[Remarks~12.4.1,
12.4.11]{monalg-rev}. Hence, using
Eqs.~\eqref{feb28-22}-\eqref{feb28-22-bis}, it follows
that $F$ contains $s$ linearly 
independent vectors in $\mathcal{A}'$, and $F$ is a facet of ${\rm RC}(I)$. 

Conversely, let $H_{(\gamma_p,-d_p)}^+$ be any of the
non-trivial closed halfspaces that occur in the irreducible
representation of ${\rm RC}(I)$ in Eq.~\eqref{supp-hyp} of
Theorem~\ref{irred-rep-rc}. Then, by \cite[Theorem 3.5]{matrof}, 
$$(\gamma_p,-d_p)=(e_{i_1}+\cdots+e_{i_k},-f_p)$$
for some $1\leq
i_1<\cdots<i_k\leq s$. Let $\mathcal{C}$ be the clutter with edge
ideal $I$ and let $\mathcal{C}^\vee$ be its blocker. Consider the induced subclutter
$\mathcal{D}=\mathcal{C}^\vee[\{t_{i_1},\ldots,t_{i_k}\}]$ of $\mathcal{C}^\vee$ consisting
of all minimal vertex covers $\{t_{j_1},\ldots,t_{j_{s-d+1}}\}$ of $\mathcal{C}$ contained in
$\{t_{i_1},\ldots,t_{i_k}\}$. Then, $I(\mathcal{D})=I_{s-d+1,k}$. Using Lemma~\ref{aug25-07} and
\cite[Theorem~2.6]{symboli}, it follows that $\gamma_p$ is an
irreducible $d_p$-cover of $\mathcal{C}$ and
$\alpha_0(\mathcal{D})$ is equal to $d_p$, where
$\alpha_0(\mathcal{D})$ is the covering number of $\mathcal{D}$, i.e.,
 the size of the smallest minimal vertex cover of $\mathcal{D}$. In particular
$\mathcal{D}\neq\emptyset$ because $d_p\geq 1$, that is, $k\geq
s-d+1$. Using the
primary decomposition of $I_{s-d+1,k}$ (cf. Eq~\eqref{feb28-22-1}), we
get 
$$\alpha_0(\mathcal{D})=k-(s-d+1)+1=k-s+d.$$
Thus, $d_p=k-s+d$
and $H_{(\gamma_p,-d_p)}^+$ occurs in the right-hand side of
Eq.~\eqref{supp-hyp-matroid}. 

Using Eq.~\eqref{supp-hyp-matroid} and Theorem~\ref{reesclu-theorem}
we obtain that Eq.~\eqref{eq1-mar1} holds. Let $u$ and $v$ be any
vertices of $\mathcal{Q}(I)$ and $\mathcal{Q}(I^\vee)$, respectively.
Then, by Eqs.~\eqref{eq1-mar1}-\eqref{eq2-mar1} and Lemma~\ref{intersection-number}, we get
$$
\langle u, v\rangle=\langle \frac{e_{i_1}+\cdots+e_{i_k}}{k-s+d},
\frac{e_{j_1}+\cdots+e_{j_{\ell}}}{\ell-d+1}\rangle=\frac{|A\cap
B|}{(k-s+d)(\ell-d+1)}\geq\frac{s}{d(s-d+1)},
$$
where $A=\{i_1,\ldots,i_k\}$ and $B=\{j_1,\ldots,j_\ell\}$, with
equality if $A=B=\{1,\ldots,s\}$. Hence
$$\min\left\{\langle u,v\rangle\mid u\in V(\mathcal{Q}(I)),\,
v\in V(\mathcal{Q}(I^\vee))\right\}={s}/{d(s-d+1)}
$$
and by Theorem~\ref{duality-ic-resurgence}, we get
$\rho_{ic}(I)=\rho_{ic}(I^\vee)=d(s-d+1)/s$. Furthermore, by
Eqs.~\eqref{eq1-mar1}-\eqref{eq2-mar1} and
\cite[Theorem~3.2]{Waldschmidt-Bocci-etal}, 
one has
$$
\widehat{\alpha}(I)=\min\{|v|\colon v\in
V(\mathcal{Q}(I^\vee))\}=\frac{s}{s-d+1}\mbox{ and }
\widehat{\alpha}(I^\vee)=\min\{|u|\colon u\in
V(\mathcal{Q}(I))\}=\frac{s}{d},
$$
where equality is attained at $(1/(s-d+1))(1,\ldots,1)$ and
$(1/d)(1,\ldots,1)$, respectively. As $I$ is normal
\cite[Proposition~2.9]{Vi4}, one has $\rho(I)=\rho_{ic}(I)$.
\end{proof}

\begin{remark}\label{numberofvertices}
If $I=I_{d,s}$ and $2\leq d<s$, then the number of vertices of $\mathcal{Q}(I)$ and
$\mathcal{Q}(I^\vee)$ are
$$|V(\mathcal{Q}(I))|=\sum_{k=s-d+1}^s\binom{s}{k}\ \mbox{ and }\ 
|V(\mathcal{Q}(I^\vee))|=\sum_{\ell=d}^s\binom{s}{\ell}.
$$
\quad This follows from Eqs.~\eqref{eq1-mar1}-\eqref{eq2-mar1}.
\end{remark}

\section{Examples}\label{examples-section}
\begin{example}\label{ic-resurgence-example1}
Let $S=\mathbb{Q}[t_1,\ldots,t_7]$ be a polynomial ring and let
$$
I=I(G)=(t_1t_3,t_1t_4,t_2t_4,t_1t_5,t_2t_5,t_3t_5,t_1t_6,
t_2t_6,t_3t_6,t_4t_6,t_2t_7,t_3t_7,t_4t_7,t_5t_7)
$$
be the edge ideal of the graph $G$ of Figure.~\ref{figure1}. This
graph is the complement of a cycle of length $7$ and is called an  
\textit{odd antihole} in the theory of perfect graphs.
\begin{figure}[ht]
\setlength{\unitlength}{.035cm}
\begin{tikzpicture}[scale=1.2,thick]
\tikzstyle{every node}=[minimum width=0pt, inner sep=2pt]
\draw (-2,0) node (0) [draw, circle, fill=gray] {};
\draw (-2.3,0) node () {$t_7$};
\draw (-1.2469796037174672,-1.56366296493606) node (1) [draw, circle, fill=gray] {};
                        \draw (-1.54,-1.56366296493606) node () {\small $t_1$};
\draw (0.44504186791262895,-1.9498558243636472) node (2) [draw, circle, fill=gray] {};
\draw (0.042504186791262895,-1.9498558243636472) node () {\small $t_2$};			
\draw (1.8019377358048383,-0.8677674782351165) node (3) [draw, circle, fill=gray] {};
\draw (2.099377358048383,-0.8677674782351165) node () { \small $t_3$};
\draw (1.8019377358048383,0.8677674782351158) node (4) [draw, circle, fill=gray] {};
\draw (2.099377358048383,0.8677674782351158) node () { \small $t_4$};
\draw (0.44504186791262895,1.9498558243636472) node (5) [draw, circle, fill=gray] {};
\draw (0.80504186791262895,1.9498558243636472) node () {\small $t_5$};
\draw (-1.2469796037174665,1.5636629649360598) node (6) [draw, circle, fill=gray] {};
\draw (-1.02469796037174665,1.7996629649360598) node () { \small $t_6$};
\draw  (0) edge (2);
\draw  (0) edge (3);
\draw  (0) edge (4);
\draw  (0) edge (5);
\draw  (1) edge (3);
\draw  (1) edge (4);
\draw  (1) edge (5);
\draw  (1) edge (6);
\draw  (2) edge (4);
\draw  (2) edge (5);
\draw  (2) edge (6);
\draw  (3) edge (5);
\draw  (3) edge (6);
\draw  (4) edge (6);
\end{tikzpicture}
\caption{Graph $G$ is the complement of a cycle of length $7$.}\label{figure1}
\end{figure}
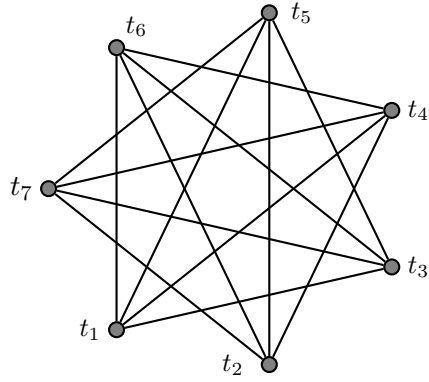

From the minimal generating set of $I$, we get that the incidence matrix $A$ of $I$ is given by
$$A={\left[{\begin{array}{cccccccccccccc}
1&1&0&1&0&0&1&0&0&0&0&0&0&0\\
0&0&1&0&1&0&0&1&0&0&1&0&0&0\\
1&0&0&0&0&1&0&0&1&0&0&1&0&0\\
0&1&1&0&0&0&0&0&0&1&0&0&1&0\\
0&0&0&1&1&1&0&0&0&0&0&0&0&1\\
0&0&0&0&0&0&1&1&1&1&0&0&0&0\\
0&0&0&0&0&0&0&0&0&0&1&1&1&1\\
\end{array}}\right]}.
$$
\quad Using Procedure~\ref{ic-resurgence-example1.m2},
\textit{Macaulay}$2$ \cite{mac2} and Theorem~\ref{duality-ic-resurgence}, we obtain the
following information. The vertex set $V(\mathcal{Q}(I))$ of $\mathcal{Q}(I)$ is given by
\begin{align*}
V(\mathcal{Q}(I))=&\{(0,\,0,\,1,\,1,\,1,\,1,\,1),\,(0,\,1,\,1,\,1,\,1,\,1,\,0),\,
(1,\,0,\,0,\,1,\,1,\,1,\,1),\,\\ 
&(1,\,1,\,0,\,0,\,1,\,1,\,1),\,
(1,\,1,\,1,\,0,\,0,\,1,\,1),\,(1,\,1,\,1,\,1,\,0,\,0,\,1),\, \\ 
&(1,\,1,\,1,\,1,\,1,\,0,\,0),\,({1}/{2},\,{1}/{2},
\,{1}/{2},\,{1}/{2},\,{1}/{2},\,{1}/{2},\,{1}/{2})\},
\end{align*}
and $\widehat{\alpha}(I^\vee)=7/2$. The incidence matrix $B$ of $I^\vee$ is
$$B=\left[\begin{matrix}
1&0&1&1&1&1&0\cr
1&1&1&1&1&0&0\cr
1&1&1&1&0&0&1\cr
1&1&1&0&0&1&1\cr
1&1&0&0&1&1&1\cr
0&1&0&1&1&1&1\cr
0&0&1&1&1&1&1
\end{matrix}
\right],
$$
the vertex set $V(\mathcal{Q}(I^\vee))$ of $\mathcal{Q}(I^\vee)$ is
given by 
\begin{align*}
V(\mathcal{Q}(I^\vee))=&\{(0,\,0,\,0,\,0,\,1,\,0,\,1),
\,(0,\,0,\,0,\,1,\,0,\,0,\,1),\,(0,\,0,\,0,\,1,\,0,\,1,\,0),\\
&\,(0,\,0,\,1,\,0,\,0,\,0,\,1),\,(0,\,0,\,1,\,0,\,0,\,1,\,0),\,(0,\,0,\,1,\,0,\,1,\,0,\,0),\\
&\,(0,\,0,\,{1}/{2},\,0,\,{1}/{2},\,0,\,{1}/{2}),
\,(0,\,1,\,0,\,0,\,0,\,0,\,1),\\
&\,(0,\,1,\,0,\,0,\,0,\,1,\,0),\,(0,\,1,\,0,\,0,\,1,\,0,\,0),\,(0,\,1,\,0,\,1,\,0,\,0,\,0),\\
&\,(0,\,{1}/{2},\,0,\,0,\,{1}/{2},\,0,\,{1}/{2}),
\,(0,\,{1}/{2},\,0,\,{1}/{2},\,0,\,0,\,{1}/{2}),\\
&\,(0,\,{1}/{2},\,0,\,{1}/{2},\,0,\,{1}/{2},\,0),\,(1,\,0,\,0,\,0,\,0,\,1,\,0),\\
&\,(1,\,0,\,0,\,0,\,1,\,0,\,0),\,(1,\,0,\,0,\,1,\,0,\,0,\,0),
\,(1,\,0,\,1,\,0,\,0,\,0,\,0),\\
&\,({1}/{2},\,0,\,0,\,{1}/{2},\,0,\,{1}/{2},\,0),
\,({1}/{2},\,0,\,{1}/{2},\,0,\,0,\,{1}/{2},\,0),\\
&\,({1}/{2},\,0,\,{1}/{2},\,0,\,{1}/{2},\,0,\,0),\,({1}/{5},\,{1}/{5},\,{1}/{5},
\,{1}/{5},\,{1}/{5},\,{1}/{5},\,{1}/{5})\},
\end{align*}
and $\widehat{\alpha}(I)=7/5$. The ic-resurgences of $I$ and
$I^\vee$ are equal to $10/7$.  
\end{example}

\begin{example}\label{example2}
Let $S=\mathbb{Q}[t_1,\ldots,t_4]$ be a polynomial ring and let
$I=({t_1}{t_2}{t_3},{t_1}{t_4},{t_2}{t_4},{t_3}{t_4})$
be the edge ideal of the clutter $\mathcal{C}$
(cf.~\cite[Example~2.26]{Francisco-TAMS}). Using
Procedure~\ref{ic-resurgence-example1.m2}, we obtain the
following information. The vertex set of
$\mathcal{Q}(I)$ is
\begin{align*}
V(\mathcal{Q}(I))=\{(0,\,0,\,1,\,1),\,(0,\,1,\,0,\,1),
\,(1,\,0,\,0,\,1),\,(1,\,1,\,1,\,0),\,({1}/{3},\,{1}/{3},\,{1}/{3},\,{2}/{3})\},
\end{align*}
$\widehat{\alpha}(I)=5/3$, $I=I^\vee$, and $\rho_{ic}(I)=9/7$.  
\end{example}

\begin{example}\label{example3}
Let $S=\mathbb{Q}[t_1,\ldots,t_7]$ be a polynomial ring and let
$$
I=(t_1t_3t_4,\, t_1t_3t_5,\, t_1t_4t_6,\, t_2t_4t_6,\, 
t_5t_6,\, t_1t_3t_7,\, t_2t_3t_7,\, t_3t_5t_7,\, t_2t_6t_7)
$$
be the edge ideal of the clutter $\mathcal{C}$. Using
Procedure~\ref{ic-resurgence-example1.m2} and Theorem~\ref{duality-ic-resurgence}, we obtain
$$
I^\vee=(t_1t_2t_5,\,  t_2t_3t_4t_5,\,  t_3t_6,\,  t_4t_5t_7,\,
t_1t_6t_7),
$$ 
$\widehat{\alpha}(I)=\widehat{\alpha}(I^\vee)=2$, and
$\rho_{ic}(I)=\rho_{ic}(I^\vee)=4/3$.  
\end{example}

\begin{example}\label{example4}
Let $S=\mathbb{Q}[t_1,\ldots,t_7]$ be a polynomial ring and let
$$
I=(t_1t_2, t_2t_3, t_3t_4, t_1t_5, t_4t_5, t_1t_6, t_2t_6, t_3t_6, t_4t_6, t_5t_6, t_1t_7, t_2t_7, t_3t_7, t_4t_7, t_5t_7, t_6t_7)
$$
be the edge ideal of the clutter $\mathcal{C}$. Using
Procedure~\ref{ic-resurgence-example1.m2} and
Theorem~\ref{duality-ic-resurgence}, we obtain
$$
I^\vee=(t_1t_2t_3t_4t_5t_6, t_1t_2t_3t_4t_5t_7, t_1t_2t_4t_6t_7,
t_1t_3t_4t_6t_7, t_1t_3t_5t_6t_7, t_2t_3t_5t_6t_7, t_2t_4t_5t_6t_7),
$$
$\alpha=({1}/{2},\,{1}/{2},\,{1}/{2},\,{1}/{2},\,{1}/{2},\,{1}/{2},\,{1}/{2})$
is a vertex of $\mathcal{Q}(I)$,
$\widehat{\alpha}(I^\vee)=|\alpha|=7/2$, 
$$\beta=({1}/{7},\,{1}/{7},\,{1}/{7},\,{1}/{7},\,{1}/{7},\,{2}/{7},\,{2}/{7})$$
is a vertex of 
$\mathcal{Q}(I^\vee)$, $\widehat{\alpha}(I)=|\beta|=9/7$, and
$\rho_{ic}(I)=\rho_{ic}(I^\vee)=14/9$.  
\end{example}

\begin{example}\label{example5}
Let $S=\mathbb{Q}[t_1,\ldots,t_6]$ be a polynomial ring and let
$I=I_{d,s}$ be the $d$-th squarefree Veronese ideal of $S$ with $d=3$
and $s=6$. The ideal $I^\vee$ is $I_{s-d+1,s}=I_{4,6}$, that is,
$I^\vee$ is generated in degree $4$,  
$\alpha=(1/4)({1},\,{1},\,{1},\,{1},\,{1},\,{1})$
is a vertex of $\mathcal{Q}(I^\vee)$,
$\widehat{\alpha}(I)=|\alpha|=3/2$, 
$$\beta=(1/3)({1},\,{1},\,{1},\,{1},\,{1},\,{1})$$
is a vertex of 
$\mathcal{Q}(I)$, $\widehat{\alpha}(I^\vee)=|\beta|=6/3=2$, and
$\rho_{ic}(I)=\rho_{ic}(I^\vee)=2$. The number of vertices of
$\mathcal{Q}(I)$ and $\mathcal{Q}(I^\vee)$ are $22$ and $42$,
respectively.  
\end{example}

\begin{appendix}

\section{Procedures}\label{Appendix}

In this appendix we give procedures for \textit{Macaulay}$2$
\cite{mac2} to compute the vertices of covering polyhedra and the
ic-resurgence of any squarefree monomial 
ideal. 

\begin{procedure}\label{ic-resurgence-example1.m2}
Let $I$ be a squarefree monomial ideal. We implement a procedure---that
uses the interface of \textit{Macaulay}$2$ \cite{mac2} 
to \textit{Normaliz} \cite{normaliz2}---to
compute the vertices of the covering polyhedra $\mathcal{Q}(I)$ and
$\mathcal{Q}(I^\vee)$, using the Rees cones of $I$ and $I^\vee$
defined in Eqs.~\eqref{rees-cone-eq} and \eqref{rees-cone-eq-dual}. Then,
using Theorem~\ref{duality-ic-resurgence}, we  
compute the ic-resurgence $\rho_{ic}(I)$ of $I$. For convenience, 
in this procedure we also include the function ``rhoichypes'' that was given in 
\cite[Algorithm~A.4]{intclos}. This procedure corresponds to
Example~\ref{ic-resurgence-example1}. To compute other examples, in
the next procedure simply change the polynomial rings $R$ and $S$, and
the generators of $I$. 
\begin{verbatim}
restart
loadPackage("Normaliz",Reload=>true)
loadPackage("Polyhedra", Reload => true)
load "SymbolicPowers.m2"
R =QQ[x1,x2,x3,x4,x5,x6,x7];
--C7-antihole
I=monomialIdeal(x1*x3,x1*x4,x1*x5,x1*x6,x2*x4,x2*x5,
x2*x6,x2*x7,x3*x5,x3*x6,x3*x7,x4*x6,x4*x7,x5*x7)
waldschmidt(I)
--Alexander dual of I
J=dual(I)
waldschmidt(J)
--transpose incidence matrix of I
A=matrix flatten apply(flatten entries gens I,
exponents)
--transpose incidence matrix of J
AJ=matrix flatten apply(flatten entries gens J,
exponents)
--generators of the Rees cone of I
M = id_(ZZ^(numcols(A)+1))^{0..numcols(A)-1}||
(A|transpose matrix {for i to numrows A-1 list 1})
--generators of the Rees cone of J
MJ=id_(ZZ^(numcols(AJ)+1))^{0..numcols(AJ)-1}||
(AJ|transpose matrix{for i to numrows AJ-1 list 1})
--rows of M
l= entries M
--rows of MJ
lJ= entries MJ
S=QQ[x1,x2,x3,x4,x5,x6,x7,x8]
L=for i in l list S_i
LJ=for i in lJ list S_i
nmzFilename="rproj1"
--Rees algebra of I
intclToricRing L
--supporting hyperplanes of the Rees cone of I
hypes=readNmzData("sup")
nmzFilename="rproj1"
--Rees algebra of J
intclToricRing LJ
--supporting hyperplanes of the Rees cone of J
hypesJ=readNmzData("sup")
--nontrivial supporting hyperplanes of RC(I)
choices = select (entries hypes, l-> 
not isSubset({last l}, {0,1}))
--nontrivial supporting hyperplanes of RC(J) 
choicesJ = select (entries hypesJ, lJ-> 
not isSubset({last lJ}, {0,1}))
A1=set choices, A2=set choicesJ
V1=apply(toList A1,toList),V2=apply(toList A2,toList)
H1=apply(V1,x->{x/-last x}),H2=apply(V2,x->{x/-last x})
--vertices of Q(I)
F1=set flatten apply(H1,
x->entries submatrix'(matrix x, ,{#x#0-1}))
--Vertices of Q(J)
F2=set flatten apply(H2,
x->entries submatrix'(matrix x, ,{#x#0-1}))
--Cartesian product Q(I) x Q(J)
E=apply(toList (F1**F2),toList)
f=(x)->matrix{x#0}*transpose(matrix{x#1})
--This is the ic-resurgence of I
rhoic=1/min toList set flatten flatten apply(apply(E,f),
entries)
--Now we compute the ic-resurgence of I using 
--the function ``rhoichypes'' 
rhoichypes = hypes -> ( 
choices = select (entries hypes, l-> not 
isSubset({last l}, {0,1}));
possibilities = for i to #choices-1 list
( 
l = choices_i;
l' = apply( drop(l,-1)|{0, last(l), 1}, a-> a_ZZ);
b  = {0};
s = select (entries hypes, l-> isSubset({last l}, 
set {0, -1}));
s' = apply(s, l -> -l |{0,0});
s' = apply(s', a-> apply(a, c-> c_ZZ));
b = b | for i to #s'-1 list 0;
v = for i to #l-2 list 0;
t = { v | {-1, 0, 1}, v | {0, 0, -1} };
t' = { v | {0, 1, 0} };
t = apply(t, a-> apply(a, c-> c_ZZ));
b = b | { 0, 0 };
b = apply(b, c-> c_ZZ);
A = matrix({l'} | s' | t);
b = transpose matrix{b};
C = matrix t';
d = matrix{{1}};
P = polyhedronFromHData( A, b, C, d);
vert = vertices P;
max flatten entries vert^{numrows vert -3});
return max possibilities
)
--This is the ic-resurgence of I
time rhoichypes hypes 
\end{verbatim}
\end{procedure}

\end{appendix}
\section*{Acknowledgments} 
We used \textit{Normaliz} \cite{normaliz2} and \textit{Macaulay}$2$
\cite{mac2} to compute the vertices of covering polyhedra
and the ic-resurgence of edge ideals. We thank the referees for a
careful reading of the paper, for pointing out
\cite[Theorem~4.8]{Geramita-Harbourne} and
\cite[Theorem~C]{Lampa-Malara}, and for the improvements suggested.


\bibliographystyle{plain}

\end{document}